\documentclass[12,reqno]{amsart}
\usepackage{amsfonts}
\usepackage{amssymb}
\usepackage{a4wide}
\usepackage{url}

\usepackage{hyperref}
\numberwithin{equation}{section}

\newtheorem{defi}{Definition}[section]
\newtheorem{thm}[defi]{Theorem}
\newtheorem{lemm}[defi]{Lemma}
\newtheorem{remark}[defi]{Remark}
\newtheorem{cor}[defi]{Corollary}

\newcommand{\R}{\mathbb{R}}

\DeclareMathOperator{\sgn}{sgn}

\begin{document}

\title[Decomposition of order statistics  of semimartingales ]
{Decomposition of order statistics of semimartingales using local times}
\author{Raouf Ghomrasni}
\address{Programme in Advanced Mathematics of Finance, School of Computational \& Applied Mathematics, University of the Witwatersrand, Private Bag 3, Wits, 2050 Johannesburg, South Africa.}
\email{Raouf.Ghomrasni@wits.ac.za}
\author{Olivier Menoukeu Pamen}
\address{Programme in Advanced Mathematics of Finance, School of Computational \& Applied Mathematics, University of the Witwatersrand, Private Bag 3, Wits, 2050 Johannesburg, South Africa.} \email{Olivier.MenoukeuPamen@students.wits.ac.za}
\thanks{{\it 2000 Mathematical Subjects Classifications.}
                               Primary 60J65, Secondary 60H05}

\date{}

\begin{abstract}

In a recent work \cite{BG}, given a collection of continuous semimartingales, authors derive a semimartingale decomposition from the corresponding ranked processes in the case that the ranked processes can meet more than two original processes at the same time. This has led to a more general decomposition of ranked processes. In this paper, we derive a more general result for semimartingales (not necessarily continuous) using a simpler approach. Furthermore, we also give a generalization of Ouknine \cite{O1,O2} and Yan's \cite{Y1} formula for local times of ranked processes.

\end{abstract}

\maketitle

\section{Introduction}

{\par Some recent developments in mathematical finance and particularly the distribution of capital in stochastic portfolio theory have led to the necessity of understanding dynamics of the $k$th-ranked amongst $n$ given stocks, at all levels $k = 1,\cdots , n$. For
example, $k = 1$ and $k = n$ correspond to the maximum and minimum process of the collection, respectively. The problem of decomposition for the maximum of $n$ semimartingales was introduced by Chitashvili and Mania in \cite{CM}. The authors showed that the maximum process can be expressed in terms of the original processes, adjusted by local times. In \cite{F}, Fernholz, defined the more general notion of ranked processes (i.e. order statistic) of $n$ continuous It\^o processes and gave the decomposition of such processes. However, the main drawback of the latter result is that, triple points do not exist, i.e., not more than two processes coincide at the same time, almost surely. Motivated by the question of extending this decomposition to triple points (and higher orders of incidence) posed by Fernholz in Problem 4.1.13 of \cite{F1}, Banner and Ghomrasni recently in \cite{BG} developed some general formulas for ranked processes of \textsl{continuous} semimartingales. In the setting of problem 4.1.13 in \cite{F1}, they showed that the ranked processes can be expressed in terms of original processes adjusted by the local times of ranked processes. The proof of those results are based on the generalization of Ouknine's formula {\cite{O1,O2,Y1}}.

 In the present paper, we give a new decomposition of order statistics of semimartingales (i.e., not necessarily continuous) in the same setting as in \cite{BG}. The obtained results are slightly different to the one in \cite{BG} in the sense that we express the order statistics of semimartingales firstly in terms of order statistic processes adjusted by their local times and secondly in terms of original processes adjusted by their local times. The proof of this result is a modified and shortened version of the proof given in \cite{BG} is based on the homogeneity property. Furthermore, we use the theory of predictable random open sets, introduced by Zheng in \cite{Z1} and the idea of the proof of Theorem 2.2 in \cite{BG} to show that 
 $$\sum_{i=1}^{n} 1_{\{X^{(i)}(t-)=0\}} \, dX^{(i)^+}(t) = \sum_{i=1}^{n} 1_{\{X_{i}(t-)=0\}} \, dX_{i}^+(t)\,, $$ where $X_i,\,\,i=1,\cdots,n$ represent the original processes and $X^{(i)}$ represent the ranked processes. As a consequence of this result, we are independently able to derive an extension of Ouknine's formula in the case of general semimartingales. The desired generalization which is essential in the demonstration of Theorem 2.3 in \cite{BG} is not used here to prove our decomposition.

{\par The paper is organized as follows. In section 2, we prove the two different decompositions of ranked processes for general semimartingales. In section 3, after showing the above equality, 

we derive a generalization of Ouknine and Yan's formula. \par}

\par}

\section{Decomposition of Ranked Semimartingales}

We begin by giving the definition of the $k$-th rank process of a family of $n$ semimartingales.
\begin{defi}
Let $X_1, \cdots, X_n$ be semimartingales. For $1 \leq k \leq n$, the {\it k-th} rank process of $X_1, \cdots, X_n$ is defined by
\begin{equation}
X^{(k)} = \max_{1	\leq i_1<\cdots<i_k \leq n} \min (X_{i_1}, \cdots, X_{i_k}),
\end{equation}
where $1 \leq i_{1}$ and $i_k \leq n$. 
\end{defi}

Note that, according to Definition 2.1, for $t\in \R^+$, 
\begin{equation}
\max_{1\leq i \leq n} X_i(t) = X^{(1)}(t) \geq X^{(2)}(t) \geq \cdots \geq X^{(n)}(t)= \min_{1\leq i \leq n} X_i(t),
\end{equation}
so that at any given time, the values of the ranked processes represent the values of the original processes arranged in descending order (i.e. the (reverse) order statistics).

{\par The following theorem shows that the ranked processes derived from semimartingales can be expressed in terms of the original processes, adjusted by local times and also in terms of ranked processes adjusted by their local times.\par}

We shall need the following definitions:
$$S_t(k)=\{i: X_i(t-)=X^{(k)}(t-)\} \text{ and } N_t(k)=\left|S_t(k)\right|\,.$$
Then $N_t(k)$ is the number of subscript $i$ such that $X_i(t)=X^{(k)}(t)$ at time $t$.
We can give a more explicit decomposition as follows:

\begin{thm}\label{theoN}
Let $X_1, \cdots, X_n$ be semimartingales. Then the {\it k-th} ranked processes $X^{(k)}, k=1,\cdots, n$, are semimartingales and we have:
\begin{align}
dX^{(k)}(t) \,=\, &  \, \sum_{i=1}^{n}  \frac{1}{N_t(k)} 1_{\{ X^{(k)}(t-)=X^{(i)}(t-)\}}\, \, dX^{(i)}(t)+\, 
\sum_{i=k+1}^{n} \frac{1}{N_t(k)}\,d\,\mathcal{L}_t^{0}((X^{(k)}-X^{(i)}))\,\nonumber \\
& - \,\sum_{i=1}^{k-1} \frac{1}{N_t(k)}\,d\,\mathcal{L}_s^0(X^{(i)}-X^{(k)}),\label{eq1cor21}\\
 \,=\, & \, \sum_{i=1}^{n}  \frac{1}{N_t(k)} 1_{\{ X^{(k)}(t-)=X_i(t-)\}}\, \, dX_i(t)+\, 
\sum_{i=1}^{n} \frac{1}{N_t(k)}\,d\,\mathcal{L}_t^{0}((X^{(k)}-X_i)^+)\,\nonumber \\
& - \,\sum_{i=1}^{n} \frac{1}{N_t(k)}\,d\,\mathcal{L}_s^0((X^{(k)}-X_i)^-),\label{eq1cor211}
\end{align}
where $\mathcal{L}_t^0 (X)=\dfrac{1}{2}L_t^0 (X)+\underset{s\leq t}{\sum}1_{\left\{X_{s-}=0\right\}}\Delta X_s$ and $L_t^0 (X)$ is the local time of the semimartingale $X$ at $0$ defined by
$$
\vert X_t\vert=\vert X_0\vert + \int_0^t \sgn(X_{s-})\, dX_s + L_t^0 (X)+\sum_{s\leq t} (\vert X_s\vert - \vert X_{s-}\vert - \sgn(X_{s-}) \Delta X_s),
$$
where $\sgn(x)=-1_{(-\infty,0]}(x)+1_{(0,\infty)}(x)$.
\end{thm}

 \proof  For all $t>0$ using the fact that we can define $N_t(k)$ as:
 $$
 N_t(k)=\sum_{i=1}^{n}1_{\left\{X^{(k)}(t-)=X_i(t-)\right\}}
=\sum_{i=1}^{n}1_{\left\{X^{(k)}(t-)=X^{(i)}(t-)\right\}}\quad a.s.\,.
 $$ 
We have the following equalities:
\begin{eqnarray}\label{eqN}
N_t(k)dX^{(k)}(t)&=&\sum_{i=1}^{n}1_{\left\{X^{(k)}(t-)=X_i(t-)\right\}}dX^{(k)}(t)\nonumber\\
&=&\sum_{i=1}^{n}1_{\left\{X^{(k)}(t-)=X^{(i)}(t-)\right\}}dX^{(k)}(t)\,.
\end{eqnarray}
By homogeneity, to show (\ref{eq1cor21}), it suffices to show that:
\begin{align}\label{eq1cor22+}
N_t(k)dX^{(k)}(t) \,=\, & \, \sum_{i=1}^{n}\,  1_{\{ X^{(k)}(t-)=X^{(i)}(t-)\}}\, \, dX^{(i)}(t)+\, 
\sum_{i=k+1}^{n} \,d\,\mathcal{L}_t^{0}((X^{(k)}-X^{(i)}))\,\nonumber \\
& - \,\sum_{i=1}^{k-1} \,d\,\mathcal{L}_s^0(X^{(i)}-X^{(k)}).
\end{align}
By the second equality of (\ref{eqN}) , we have 
\begin{eqnarray*}
N_t(k)dX^{(k)}(t)&=&\sum_{i=1}^{n}1_{\left\{X^{(k)}(t-)=X^{(i)}(t-)\right\}}dX^{(k)}(t)\\
&=&\sum_{i=1}^{n}1_{\left\{X^{(k)}(t-)=X^{(i)}(t-)\right\}}dX^{(i)}(t)\\
&&+\sum_{i=1}^{n}1_{\left\{X^{(k)}(t-)=X^{(i)}(t-)\right\}}d\left(X^{(k)}(t)-X^{(i)}(t)\right)\,.
\end{eqnarray*}
Using the formula
\begin{eqnarray}\label{eqlocalt1}
\mathcal{L}_t^{0}(Z)=\int_0^t\,1_{\left\{Z(s-)=0\right\}}dZ(s)\,,
\end{eqnarray}

which is valid for nonnegative semimartingales $Z$.\\

By applying (\ref{eqlocalt1}) to $N_t(k)dX^{(k)}(t),\,\,t>0$, we obtain:

\begin{align}\label{eqlocalt2}
N_t(k)dX^{(k)}(t)\,=\,&\sum_{i=1}^{n}1_{\left\{X^{(k)}(t-)=X^{(i)}(t-)\right\}}dX^{(i)}(t)\nonumber\\
&+\,\sum_{i=1}^{n}1_{\left\{X^{(k)}(t-)=X^{(i)}(t-)\right\}}d\left((X^{(k)}(t)-X^{(i)}(t))^+\right)\nonumber\\
&-\,\sum_{i=1}^{n}1_{\left\{X^{(k)}(t-)=X^{(i)}(t-)\right\}}d\left((X^{(k)}(t)-X^{(i)}(t))^-\right)\nonumber\\
\,=\,& \sum_{i=1}^{n}1_{\left\{X^{(k)}(t-)=X^{(i)}(t-)\right\}}dX^{(i)}(t)\,+\,\sum_{i=1}^{n}\,d\,\mathcal{L}_t^0\left((X^{(k)}-X^{(i)})^+\right)\nonumber\\
& -\,\sum_{i=1}^{n}\,d\,\mathcal{L}_t^0\left((X^{(k)}-X^{(i)})^-\right)\,.
\end{align}
Nothing that:
\begin{equation*}
(X^{(k)}-X^{(j)})^+=
\left\{\begin{array}{lll}
X^{(k)}-X^{(j)},&\text{if }j>k\\
0,&\text{if }j\leq k
\end{array}
\right.
\end{equation*}
and that
\begin{equation*}
(X^{(k)}-X^{(j)})^-=
\left\{\begin{array}{lll}
X^{(j)}-X^{(k)},&\text{if }j<k\\
0,&\text{if }j\geq k
\end{array}
\right.
\end{equation*}
then (\ref{eq1cor21}) follows.

In the same way, we prove (\ref{eq1cor211}) by applying the first equality of (\ref{eqN}), and (\ref{eqlocalt1}).

\qed

\subsection{Local time and Norms}

\text{ }\\
\text{ }

The next result is proved in \cite{CO}:
\begin{lemm}
Let $X=X_1,\cdots,X_n$ be a $n-$dimensional semimartingale, $N_1$ and $N_2$ be norms on $\mathbb{R}^n$ such that $N_1\leq N_2$. Then $L^0_t(N_1(X))=L^0_t(N_2(X))$. 
\end{lemm}
For example

$$
L^0_t(\underset{1\leq i\leq n}{\max}\left|X_i\right|)\leq L^0_t(\sum_{i=1}^{n}\left|X_i\right|)\leq n\,L^0_t(\underset{1\leq i\leq n}{\max}\left|X_i\right|).
$$

For positive continuous semimartingales, we have the following result.
\begin{cor}
Let $X_1,\cdots,X_n$ be positive continuous semimartingales. The we have inequality hold:
$$
L^0_t(\sum_{i=1}^{n}X_i)\leq n\,\sum_{i=1}^{n}L^0_t(X_i).
$$
\end{cor}
\proof
It is known that, the next equality holds for continuous semimartigales (see \cite{BG}): $\sum_{i=1}^{n} L_t^0(X^{(i)}) = \sum_{i=1}^{n} L_t^0(X_{i}) \,\, (\ast)$, putting $L_t^0(X^{(1)})=L^0_t(\underset{1\leq i\leq n}{\max}X_i)$, we have by the preceding lemma 
\begin{eqnarray*}
L^0_t(\sum_{i=1}^{n}X_i)&\leq & n\,L^0_t(\underset{1\leq i\leq n}{\max}X_i)=n\,L_t^0(X^{(1)})\\
&\leq& n\, \sum_{i=1}^{n} L_t^0(X_{i}) \,\,\text{ }\,(\text{ by } (\ast))\,.
\end{eqnarray*}

\qed

\begin{remark}
In theorem \ref{theoN}, consider the first rank process $X^{(1)}$ i.e., the maximum process, then for a continuous semimartingale, the equalities (\ref{eq1cor21}) and (\ref{eq1cor211}) become: 
\begin{align}
dX^{(1)}(t) \,=\, &  \, \sum_{i=1}^{n}  \frac{1}{N_t(1)} 1_{\{ X^{(1)}(t)=X^{(i)}(t)\}}\, \, dX^{(i)}(t)+\, 
\sum_{i=2}^{n} \frac{1}{N_t(1)}\,d\,L_t^{0}((X^{(1)}-X^{(i)}))\,\nonumber \\
 \,=\, & \, \sum_{i=1}^{n}  \frac{1}{N_t(1)} 1_{\{ X^{(1)}(t)=X_i(t)\}}\, \, dX_i(t)+\, 
\sum_{i=1}^{n} \frac{1}{N_t(1)}\,d\,L_t^{0}((X^{(1)}-X_i))\,,\nonumber
\end{align}
and by homogeneity we have 
\begin{align}
N_t(1)dX^{(1)}(t) \,=\, &  \, \sum_{i=1}^{n}   1_{\{ X^{(1)}(t)=X^{(i)}(t)\}}\, \, dX^{(i)}(t)+\, 
\sum_{i=2}^{n} \,d\,L_t^{0}((X^{(1)}-X^{(i)}))\,\nonumber \\
 \,=\, & \, \sum_{i=1}^{n}  1_{\{ X^{(1)}(t)=X_i(t)\}}\, \, dX_i(t)+\, 
\sum_{i=1}^{n} \,d\,L_t^{0}((X^{(1)}-X_i))\,.\label{eqrem1}
\end{align}
Define the processes $Y_1, \cdots, Y_n$ by: $Y_i(t)=X^{(1)}(t)-X_i(t), \,\,i=1,\cdots,n$ then $Y_1, \cdots, Y_n$ are continuous semimartingales and the processes $Y^{(1)}, \cdots, Y^{(n)}$ defined by $Y^{(i)}(t)=X^{(1)}(t)-X^{(i)}(t), \,\,i=1,\cdots,n$ are the {\it i-th} ranked processes of $Y_i(t), \,\,i=1,\cdots,n$ with the property $Y^{(1)}\leq Y^{(2)}\leq \cdots \leq Y^{(n)}$ and there are continuous semimatingales. It is shown by Banner and Ghomrasni in \cite{BG} that
\begin{equation*}
\sum_{i=1}^{n} L_t^0(Y^{(i)}) = \sum_{i=1}^{n} L_t^0(Y_{i})\,,
\end{equation*} which is a generalization of Ouknine and Yan's formula for continuous semimartingles.
Replacing $Y^{(i)}$ and $Y_{i}$ by their expressions, we have:
  $$
\sum_{i=1}^{n} L_t^0(X^{(1)}(t)-X^{(i)}) = \sum_{i=1}^{n} L_t^0(X^{(1)}(t)-X_i(t))\,.
$$
Using this equation and identifying the first and the second equalities of (\ref{eqrem1}), we conclude that:
$$
\sum_{i=1}^{n}   1_{\{ X^{(1)}(t)=X^{(i)}(t)\}}\, \, dX^{(i)}(t)\,=\,\sum_{i=1}^{n}  1_{\{ X^{(1)}(t)=X_i(t)\}}\, \, dX_i(t)\,.
$$
\end{remark}
The questions are: Does such an equality hold if we replace $X^{(1)}$ by $X^{(k)}$ for an arbitrary $k\in \{2,\cdots,n\}$? Is there any type of this equality for a general semimartingale? The answers of theses questions are given in the next section.

\section{Generalization of Ouknine and Yan's Formula for Semimartingales}

In this section we derive a generalization of Ouknine and Yan's formula for semimartingales. Such a result was proved in \cite{BG} in the case of continuous semimartingales. In order to give such an extension, we need first to prove the next theorem.

\begin{thm}\label{theosum}
 Let $X_1, \cdots, X_n$ be semimartingales. Then the following equality hold: 
\begin{equation}
\sum_{i=1}^{n} 1_{\{X^{(i)}(t-)=0\}} \, dX^{(i)^{+}}(t) = \sum_{i=1}^{n} 1_{\{X_{i}(t-)=0\}} \, dX^{+}_{i}(t).
\end{equation}
\end{thm}

\proof  We will proceed by induction. The case $n=1$ is trivial. For $n=2$, let show that 
\begin{eqnarray}\label{eqsum2}
&&1_{\{X^{(1)}(t-)=0\}} \, dX^{(1)^+}(t)\,+\,1_{\{X^{(2)}(t-)=0\}} \, dX^{(2)^+}(t)\notag\\
&=&\, 1_{\{X_{1}(t-)=0\}} \, dX_{1}^+(t)\,+\,1_{\{X_{2}(t-)=0\}} \, dX_{2}^+(t)\,,
\end{eqnarray}
where $X^{(1)}=X_1 \vee X_2$ and $X^{(2)}=X_1 \wedge X_2$. At this point we follows the same idea as in the proof of the second theorem in \cite{O1}. Since 
\begin{align}
\left\{X_1 (t-)\vee X_2(t-)=0\right\}\,=\,&\left\{X_1(t-)<X_2(t-)=0\right\}\cup \left\{X_2(t-)<X_1(t-)=0\right\}\notag \\
\,&\cup\left\{X_1(t-)=X_2(t-)=0\right\}\,,\notag
\end{align} and 
\begin{align}
\left\{X_1(t-) \wedge X_2(t-)=0\right\}\,=\,&\left\{X_1(t-)>X_2(t-)=0\right\}\cup \left\{X_2(t-)>X_1(t-)=0\right\}\notag\\
&\cup\left\{X_1(t-)=X_2(t-)=0\right\}\,.\notag
\end{align}
 We can write:
\begin{eqnarray}\label{eqmax2}
&&1_{\{X^{(1)}(t-)=0\}} \, dX^{(1)^+}(t)\notag\\
&=&1_{\left\{X_1(t-)<X_2(t^-)=0\right\}}dX^{(1)^{+}}(t)\,+\,1_{\left\{X_2(t-)<X_1(t-)=0\right\}}dX^{(1)^{+}}(t)\notag\\
&&+\,1_{\left\{X_1(t-)=X_2(t-)=0\right\}}dX^{(1)^{+}}(t)\,,
\end{eqnarray}
and 
\begin{eqnarray}\label{eqmin2}
&&1_{\{X^{(2)}(t-)=0\}} \, dX^{(2)^{+}}(t)\notag\\
&=&1_{\left\{X_1(t-)>X_2(t^-)=0\right\}}dX^{(2)^{+}}(t)\,+\,1_{\left\{X_2(t-)>X_1(t-)=0\right\}}dX^{(2)^{+}}(t)\notag\\
&&+\,1_{\left\{X_1(t-)=X_2(t-)=0\right\}}dX^{(2)^{+}}(t)\,.
\end{eqnarray}

In the predictable random open set $\left\{X_1<X_2\right\}$ the semimartingales $X^{(1)}$ and $X_2^{+}$ are equal, then the first term of the right hand side of (\ref{eqmax2}) is $1_{\left\{X_1(t-)<0\right\}}1_{\left\{X_2(t-)=0\right\}}dX^{+}_2(t)$. Applying the same reasoning, the second term is $1_{\left\{X_2(t-)<0\right\}}1_{\left\{X_1(t-)=0\right\}}dX^{+}_1(t)$. For the third term, we can write $X^{(1)^{+}}=(X_1 \vee X_2)^{+}=X_1^+ \vee X_2^+=X_2^++(X_1^+-X_2^+)^+$ and then it becomes $$1_{\left\{X_1(t-)=X_2(t-)=0\right\}}dX_2^+(t)\,+\,1_{\left\{X_1(t-)=X_2(t-)=0\right\}}d(X_1^+(t)-X_2^+(t))^+ \,. $$
These remarks allow us to write (\ref{eqmax2}) as:
\begin{eqnarray}\label{eqmax21}
&&1_{\{X^{(1)}(t-)=0\}} \, dX^{(1)^+}(t) \notag\\
&=& 1_{\left\{X_1(t-)<0\right\}}1_{\left\{X_2(t^-)=0\right\}}dX^+_2(t)\,+\,1_{\left\{X_2(t^-)<0\right\}}1_{\left\{X_1(t^-)=0\right\}}dX^+_1(t)\notag\\
&&\,+\,1_{\left\{X_1(t^-)=X_2(t)=0\right\}}dX^+_2\,+\,1_{\left\{X_1(t^-)=X_2(t^-)=0\right\}}d(X_1^+-X_2^+)^+\notag\\
&=&1_{\left\{X_1(t-)<0\right\}}1_{\left\{X_2(t^-)=0\right\}}dX_2^+(t)\,+\,1_{\left\{X_2(t-)<0\right\}}1_{\left\{X_1(t-)=0\right\}}dX_1^+(t)\notag\\
&&\,+\,1_{\left\{X_1(t-)=0\right\}}1_{\left\{X_2(t-)=0\right\}}dX_2^+(t)\,+\,1_{\left\{X_1(t-)=X_2(t-)=0\right\}}d(X_1^+(t)-X_2(t)^+)^+\,,
\end{eqnarray}

Following the argument for the process $X^{(2)}$, we obtain:
\begin{eqnarray}\label{eqmin21}
&&1_{\{X^{(2)}(t-)=0\}} \, dX^{(2)^+}(t)\notag\\
&=& 1_{\left\{X_1(t-)>0\right\}}1_{\left\{X_2(t-)=0\right\}}dX_2^+(t)\,+\,1_{\left\{X_2(t-)>0\right\}}1_{\left\{X_1(t-)=0\right\}}dX_1^+(t)\notag \\
&&\,+\,1_{\left\{X_1(t-)=0\right\}}1_{\left\{X_2(t-)=0\right\}}dX_1^+(t)\,-\,1_{\left\{X_1(t-)=X_2(t-)=0\right\}}d(X_1^+(t)-X_2^+(t))^-\,,
\end{eqnarray}

where we have used the fact that $X^{(1)^+}=(X_1 \wedge X_2)^+=X_1^+ \wedge X_2^+=X_1^+-(X_1^+-X_2^+)^+$.\\

Summing (\ref{eqmax21}) and (\ref{eqmin21}) we obtain the desired result for $n=2$.\\
Now assume the result holds for some $n$. We adjust here the proof given by Banner and Ghomrasni in \cite{BG}. Given semimartingales $X_1,\cdots,X_n,X_{n+1}$, we define $X^{(k)},\,k=1,\cdots,n$, as above and also set
$$
X^{[k]}(\cdot) = \max_{1\leq i_1<\cdots<i_k\leq n+1} \min (X_{i_1}(\cdot), \cdots, X_{i_k}(\cdot))\,.
$$

The process $X^{[k]}(\cdot)$ is the $k$th-ranked process with respect to all $n+1$ semimartingales $X_1,\cdots,X_n,X_{n+1}$. It will be convenient to set $X^{(0)}(\cdot):\equiv \infty$. In order to show the equality for $n+1$ we are starting by showing that:

\begin{eqnarray}\label{eqiden}
1_{\left\{X^{(k-1)}(t-)\wedge X_{n+1}(t-)=0\right\}}d(X^{(k-1)^+}(t)\wedge X_{n+1}^+(t))+1_{\left\{X^{(k)}(t-)=0\right\}}dX^{(k)^+}(t)\notag\\
= 1_{\left\{X^{[k]}(t-)=0\right\}}dX^{[k]^+}(t)
+1_{\left\{X^{(k)}(t-)\wedge X_{n+1}(t-)=0\right\}}d(X^{(k)^+}(t)\wedge X_{n+1}^+(t))
\end{eqnarray}

for $k=1,\cdots,n$ and $t>0$. Suppose first that $k>1$. By (\ref{eqsum2}), we have 
 \begin{eqnarray*}
1_{\left\{X^{(k-1)}(t-)\wedge X_{n+1}(t-)=0\right\}}d(X^{(k-1)^+}(t)\wedge X_{n+1}^+(t))+1_{\left\{X^{(k)}(t-)=0\right\}}dX^{(k)^+}(t)\\=
1_{\left\{(X^{(k-1)}(t-)\wedge X_{n+1}(t^-))\vee X^{(k)}(t^-)=0\right\}}d\left((X^{(k-1)^+}(t)\wedge X_{n+1}^+(t))\vee X^{(k)^+}(t)\right)\\
+1_{\left\{(X^{(k-1)}(t-)\wedge X_{n+1}(t-))\wedge X^{(k)}(t^-)=0\right\}}d\left((X^{(k-1)^+}(t)\wedge X_{n+1}^+(t))\wedge X^{(k)^+}(t)\right)\,.
 \end{eqnarray*}
 Since $X^{(k)}(t) \leq X^{(k-1)}(t)$ for all $t>0$, the second term of the right hand side of the above equation is simply 
 $1_{\left\{ X_{n+1}(t^-)\wedge X^{(k)}(t-)=0\right\}}d( X_{n+1}^+(t)\wedge X^{(k)^+}(t))$. On the other hand, we have 
 \begin{equation*}
 (X^{(k-1)}\wedge X_{n+1})\vee X^{(k)}(t)=\left\{
 \begin{array}{lll}
 X^{(k-1)}(t) &\text{ if } X_{n+1}(t)\geq X^{(k-1)}(t)\geq X^{(k)}(t)\\
  X_{n+1}(t)&\text{ if } X^{(k-1)}(t)\geq X_{n+1}(t) \geq X^{(k)}(t) \\
 X^{(k)}(t) &\text{ if }  X^{(k-1)}(t) \geq X^{(k)}(t) \geq X_{n+1}(t)
 \end{array}
 \right. 
 \end{equation*}
 
 In each case it can be checked that $(X^{(k-1)}\wedge X_{n+1})\vee X^{(k)}(t)$ is the $k$th smallest of the numbers $X_1,\cdots,X_{n+1}$; that is, $(X^{(k-1)}\wedge X_{n+1})\vee X^{(k)}(\cdot)\equiv X^{[k]}(\cdot)$. It follows that $X^{[k]}$ is a continuous semimartingale for $k=1,\cdots,n$. Equation (\ref{eqiden}) follows for $k=2,\cdots,n$. If $k=1$, then $X^{(0)}(\cdot)\equiv \infty$, applying (\ref{eqsum2}), (\ref{eqiden}) reduces to
 \begin{eqnarray}
&&1_{\left\{ X_{n+1}(t^-)=0\right\}}d X_{n+1}^+(t)+1_{\left\{X^{(1)}(t-)=0\right\}}dX^{(1)^+}(t)\notag \\
\,&=&\,1_{\left\{X^{(1)}(t-)\vee X_{n+1}(t-)=0\right\}}d(X^{(1)^+}(t)\vee X_{n+1}^+(t))\notag\\
\,&+&\,1_{\left\{X^{(1)}(t-)\wedge X_{n+1}(t-)=0\right\}}d(X^{(1)^+}(t)\wedge X_{n+1}^+(t))\notag\\
\,&=&\,1_{\left\{X^{(1)}(t-)\wedge X_{n+1}(t-)=0\right\}}d(X^{(1)^+}(t)\wedge X_{n+1}^+(t))\notag\\
\,&+&\,1_{\left\{ X^{[1]}(t-)=0\right\}}d X^{[1]^+}(t)\,,
\end{eqnarray}

where we observed that $(X^{(1)}\vee X_{n+1})(\cdot)\equiv X^{[1]}(\cdot)$.\\
 Finally, by the induction hypothesis and (\ref{eqiden}) , we have
 \begin{eqnarray*}
\sum_{i=1}^{n+1} 1_{\{X_{i}(t-)=0\}} \, dX_{i}^+(t)&=&\sum_{i=1}^{n} 1_{\{X_{i}(t-)=0\}} \, dX_{i}^+(t)\,+\,1_{\{X_{n+1}(t-)=0\}} \, dX_{n+1}^+(t)\\
&=&\sum_{i=1}^{n} 1_{\{X^{(i)}(t-)=0\}} \, dX^{(i)^+}(t)\,+\,1_{\{X_{n+1}(t-)=0\}} \, dX_{n+1}^+(t)\\
&=&\sum_{i=1}^{n} 1_{\{X^{[i]}(t-)=0\}} \, dX^{[i]^+}(t)\,+\,1_{\{X_{n+1}(t-)=0\}} \, dX_{n+1}^+(t)\\
&&-\,\sum_{i=1}^{n}1_{\left\{X^{(i-1)}(t-)\wedge X_{n+1}(t-)=0\right\}}d(X^{(i-1)^+}(t)\wedge X_{n+1}^+(t))\\
&&+\,\sum_{i=1}^{n}1_{\left\{X^{(i)}(t-)\wedge X_{n+1}(t-)=0\right\}}d(X^{(i)^+}(t)\wedge X_{n+1}^+(t))\\
&=&\sum_{i=1}^{n} 1_{\{X^{[i]}(t-)=0\}} \, dX^{[i]^+}(t)\,+\,1_{\{X_{n+1}(t-)=0\}} \, dX_{n+1}^+(t)\\
&&-1_{\left\{X^{(0)}(t-)\wedge X_{n+1}(t-)=0\right\}}d(X^{(0)^+}(t)\wedge X_{n+1}^+(t))\\
&&+\,1_{\left\{X^{(n)}(t-)\wedge X_{n+1}(t-)=0\right\}}d(X^{(i)^+}(t)\wedge X_{n+1}^+(t))\\
&=&\sum_{i=1}^{n+1} 1_{\{X^{[i]}(t-)=0\}} \, dX^{[i]^+}(t)\,.
 \end{eqnarray*}
 The third equality follows from (\ref{eqiden}) while the last come from the fact that \\$X^{(0)}(t)\wedge X_{n+1}(t)= X_{n+1}(t)$ and $(X^{(n)}\wedge X_{n+1})(\cdot)\equiv X^{[n+1]}(\cdot)$ for all $t>0$; then the result follows by induction
\qed

In the case of continuous semimartingales, the preceding theorem becomes:
\begin{cor}\label{corsum}
 Let $X_1, \cdots, X_n$ be continuous positive semimartingales. Then the following equality holds: 
\begin{equation}
\sum_{i=1}^{n} 1_{\{X^{(i)}(t)=0\}} \, dX^{(i)}(t) = \sum_{i=1}^{n} 1_{\{X_{i}(t)=0\}} \, dX_{i}(t).
\end{equation}
\end{cor}
It follows that:

\begin{cor}\label{corollary1}
Let $X_1, \cdots, X_n$ be continuous semimartingales. Then the {\it k-th} ranked processes $X^{(k)}, k=1,\cdots, n$, are continuous semimartingales and we have:
\begin{eqnarray}\label{eqlem1}
\sum_{i=1}^{n}1_{\left\{X^{(k)}(t)=X^{(i)}(t)\right\}}d\Big(X^{(k)}(t)-X^{(i)}(t)\Big)^+\,=\,\sum_{i=1}^{n}1_{\left\{X^{(k)}(t)=X_i(t)\right\}}d\Big(X^{(k)}(t)-X_i(t)\Big)^+\,.
\end{eqnarray}
\end{cor}

\proof  Fix $X^{(k)}$, for $k=1,\cdots,n$ and define the processes $Y_1, \cdots, Y_n$ by: $Y_i(t)=X^{(k)}(t)-X_i(t), \,\,i=1,\cdots,n$ then $Y_1, \cdots, Y_n$ are continuous semimartingales and the processes $Y^{(1)}, \cdots, Y^{(n)}$ defined by $Y^{(i)}(t)=X^{(k)}(t)-X^{(i)}(t), \,\,i=1,\cdots,n$ are the {\it i-th} ranked processes of $Y_i(t), \,\,i=1,\cdots,n$ with the property $Y^{(1)}\leq Y^{(2)}\leq \cdots \leq Y^{(n)}$ and there are continuous semimatingales. 
By theorem \ref{theosum}, we have 
\begin{equation}
\sum_{i=1}^{n} 1_{\{Y^{(i)}(t)=0\}} \, dY^{(i)}(t) = \sum_{i=1}^{n} 1_{\{Y_{i}(t)=0\}} \, dY_{i}(t) \,,
\end{equation}
and the result follows
\qed

A consequence of Theorem \ref{theosum} is the following theorem, which is a generalization of Yan \cite{Y1}, Ouknine's \cite{O1, O2} formula.

\begin{thm}\label{theorem2}
Let $X_1, \cdots, X_n$ be semimartingales. Then we have:
\begin{equation}
\sum_{i=1}^{n} L_t^0(X^{(i)}) = \sum_{i=1}^{n} L_t^0(X_{i}),
\end{equation}
where $L_t^0 (X)$ is the local time of the continuous semimartingale $X$ at $0$.
\end{thm}

\proof We recall first that $L_t^0(Z)=L_t^0(Z^+)$ for every semimartingale $Z$. By theorem \ref{theosum} the following equality holds:
\begin{equation}
\sum_{i=1}^{n} 1_{\{X^{(i)}(t-)=0\}} \, dX^{(i)^+}(t) = \sum_{i=1}^{n} 1_{\{X_{i}(t-)=0\}} \, dX_{i}^+(t)\,.
\end{equation}
We know that 
$$
\mathcal{L}_t^{0}(Z^+)=\int_0^t\,1_{\left\{Z(s-)=0\right\}}dZ^+(s)\,,
$$ 
for all semimartingales. Then the preceding equation becomes
\begin{equation}
\sum_{i=1}^{n} \mathcal{L}_t^{0}(X^{(i)^+}) = \sum_{i=1}^{n} \mathcal{L}_t^{0}( X_{i}^+)\,.
\end{equation}
Putting 
\begin{eqnarray*}
A(t)=\sum_{i=1}^{n} \mathcal{L}_t^{0}(X^{(i)^+}) \text{ and } B(t)=\sum_{i=1}^{n} \mathcal{L}_t^{0}(X_i^+)\,.\\
\end{eqnarray*}
 then
  \begin{eqnarray*}
 A(t)&=&\sum_{i=1}^{n} \left(\dfrac{1}{2}L_t^0 (X^{(i)^+})+\underset{s\leq t}{\sum}1_{\left\{X^{(i)}(s-)=0\right\}}\Delta X^{(i)^+}(s)\right)\,,\\
 B(t)&=&\sum_{i=1}^{n} \left(\dfrac{1}{2}L_t^0 (X_{i}^+)+\underset{s\leq t}{\sum}1_{\left\{X_{i}(s-)=0\right\}}\Delta X_{i}^+(s)\right)\,.
 \end{eqnarray*}
 Since $A(t)=B(t)$ for all $t>0$, we have $A^c(t)=B^c(t)$ where $A^c$ (resp. $B^c$) is the continuous part of $A$ (resp. $B$). The desired result follows from the continuity of local time and the fact $L_t^0(Z)=L_t^0(Z^+)$.
\qed

In particular,

\begin{cor}{\rm \bf Yan \cite{Y1}, Ouknine \cite{O1, O2}}
\\
Let $X$ and $Y$ be semimartingales. It is shown that 
\begin{equation}
L^0_t(X\vee Y)+ L^0_t(X\wedge Y)=L^0_t(X)+L^0_t(Y),
\end{equation}
where $L^0_t(X)$ $(t\ge 0)$ denotes the local time at 0 of $X$.
\end{cor}

%


\begin{thebibliography}{99}


\bibitem{BG} A.D. Banner {\it and} R. Ghomrasni, \textit{Local times of ranked continuous semimartingales}.  Stochastic Process. Appl.  {\bf 118} (2008), 1244–-1253.  

\bibitem{CM} R. J. Chitashvili {\it and} M. G. Mania, \textit{Decomposition of the maximum of semimartingales and generalized It\^o's formula}. New trends in probability and statistics, Vol. {\bf 1} (Bakuriani, 1990), 301--350.

\bibitem{CO} F. Coquet  {\it and} Y. Ouknine, \textit{Some identities on semimartingales local times}. Statist. Probab. Lett.  {\bf 49} (2000), 149--153.


\bibitem{F} R. Fernholz, \textit{Equity portfolios generated by functions of ranked market weights.} Finance Stoch. {\bf 5} (2001), 469--486.

\bibitem{F1} R. Fernholz, \textit{Stochastic portfolio theory.} Springer-Verlag, New York, 2002.


\bibitem{L} Jun S. Liu, \textit{Siegel's formula via Stein's identities.} Statist. Probab. Lett. {\bf 21} (1994), 247--251.

\bibitem{O1} Y. Ouknine, \textit{Temps local du produit et du sup de deux semimartingales.} {S\'eminaire de Probabilit\'es, XXIV, 1988/89} {Lecture Notes in Math., } {\bf 1426} (1990), 477--479.

\bibitem{O2} Y. Ouknine, \textit{G\'en\'eralisation d'un lemme de S. Nakao et applications.} Stochastics {\bf 23}
 (1988), 149--157.
 
 
\bibitem{RY} D. Revuz {\it and} M. Yor, \textit{Continuous Martingales And Brownian Motion.} Third edition, Springer-Verlag, Berlin, 2004.


\bibitem{Y1} J. A. Yan, \textit{A formula for local times of semimartingales}. Northeast. Math. J. {\bf 1} (1985), 138--140.

\bibitem{Z1} W.A. Zhen, \textit{Semimartingales in predictable random open sets}. S\'eminaire de Probabilit\'e XVI, Lecture Notes in Math., {\bf 920} (1982), 370--379.


\end{thebibliography}
\end{document}